\documentclass[reqno]{amsart}
\usepackage{amsmath, amsthm, amscd, amssymb, amsfonts, amsbsy}
\usepackage{mathrsfs}
\usepackage{pxfonts}
\usepackage{bbm}
\usepackage{todonotes}
\usepackage{marginnote}

\theoremstyle{plain}
\newtheorem{theorem}[equation]{Theorem}
\newtheorem{lemma}[equation]{Lemma}

\theoremstyle{definition}

\theoremstyle{remark}
\newtheorem{remark}[equation]{Remark}

\newcommand{\capacity}{\operatorname{cap}}

\newcommand{\supp}{\operatorname{supp}}

\numberwithin{equation}{section}

\newcommand{\bR}{\mathbb{R}}

\providecommand{\abs}[1]{\lvert#1\rvert}

\providecommand{\norm}[1]{\lVert#1\rVert}

\newcommand{\red}{\mathscr R}
\newcommand{\bal}{\hat{\mathscr R}}

\begin{document}

\title[Quantitative Wiener criterion]
{A quantitative Wiener criterion for the heat equation}

\author[S. Kim]{Seick Kim}
\address[S. Kim]{Department of Mathematics, Yonsei University, 50 Yonsei-Ro, Seodaemun-gu, Seoul 03722, Republic of Korea}
\email{kimseick@yonsei.ac.kr}
\thanks{S. Kim is supported by the National Research Foundation of Korea (NRF) under agreement NRF-2022R1A2C1003322.
S. Kim is also affiliated with the Intelligent Computational Science Institute (IN2CSI) at Yonsei University.}

\date{}


\begin{abstract}
We give a short proof of the Wiener criterion of Evans and Gariepy for the heat equation.
The proof also yields an explicit modulus of continuity at every regular boundary point in terms of the partial sums of the classical Wiener series and the modulus of continuity of the boundary data.
\end{abstract}

\maketitle

\section{Introduction and main result}
We give a short proof of the Wiener criterion of Evans and Gariepy \cite{EG} for the heat equation, supplemented by a quantitative boundary estimate.
Lanconelli \cite{Lanc75} obtained quantitative boundary estimates for a more general class of parabolic operators, but the Wiener-type series in his estimate involves a loss in the Gaussian exponent.
Consequently, even for the heat equation, his estimate yields a modulus of continuity only under a condition stronger than boundary regularity.
In contrast, our estimate is expressed in terms of the classical Wiener series of Evans and Gariepy and yields a modulus of continuity at every regular boundary point.
To the best of our knowledge, no such estimate has previously been established.
Throughout the paper, boundary regularity is understood in the Perron sense.

Let $P=\partial_t-\Delta$ and let $\Gamma$ be its fundamental solution given by
\[
\Gamma(t,x)=(4\pi t)^{-d/2} \exp\left(-\frac{\abs{x}^2}{4t}\right)\qquad(t>0),
\]
with $\Gamma(t,x)=0$ for $t\le0$.

For $z=(t,x)\in\bR^{d+1}$, we write $\Gamma(z)=\Gamma(t,x)$ and set
\[
\abs{z}=\max\left\{\abs{x},\sqrt{\abs{t}}\right\},
\]
so that $\abs{z-\zeta}$ denotes the parabolic distance between $z$ and $\zeta$.

For a compact set $K\subset\bR^{d+1}$, define its thermal capacity by
\[
\capacity(K)=\sup\{\mu(K):\mu\ge0,\ \supp\mu\subset K,\; \Gamma*\mu \le 1\},
\]
where $\mu$ ranges over Radon measures and the inequality $\Gamma*\mu \le 1$ is understood pointwise in $\bR^{d+1}$.
For an arbitrary set $E\subset\bR^{d+1}$, we define
\[
\capacity(E)=\sup_{K\Subset E}\capacity(K),
\]
where $K\Subset E$ means that $K$ is a compact subset of $E$.

For $z_0\in\bR^{d+1}$ and $r>0$, define the heat ball
\[
\mathcal E(z_0,r)=\left\{\zeta\in\bR^{d+1}:\Gamma(z_0-\zeta)>(4\pi r)^{-d/2}\right\}
\]
and the annular set
\[
\mathcal A(z_0,r)=\overline{\mathcal E(z_0,r)}\setminus\mathcal E(z_0,r/2).
\]

For $s\in\bR$, let $s_+=\max\{s,0\}$, and let $\lfloor s\rfloor$ and $\lceil s\rceil$ denote its floor and ceiling, respectively.
\begin{theorem}	\label{thm01}
Let $\Omega\subset\bR^{d+1}$ be a bounded open set and $z_0\in\partial\Omega$.
For integers $n\ge1$, set
\[
\mathscr W_n=\sum_{j=1}^n 2^{jd/2} \capacity\left(\Omega^\complement\cap\mathcal A(z_0,2^{-j})\right),
\]
and put $\mathscr W_n=0$ for $n \le 0$.
Then $z_0$ is regular for $P$ in $\Omega$ if and only if
\begin{equation}	\label{eq1501wed}
\lim_{n\to\infty}\mathscr W_n=\sum_{j=1}^\infty 2^{jd/2}\capacity\left(\Omega^\complement\cap\mathcal A(z_0,2^{-j})\right)=\infty.
\end{equation}

Moreover, for $\varphi\in C(\partial\Omega)$, let $H_\varphi^\Omega$ be the Perron solution with boundary
data $\varphi$ and set
\[
\omega_\varphi(z_0,r)=\sup_{\substack{\zeta\in\partial\Omega\\\abs{\zeta-z_0}\le r}}\abs{\varphi(\zeta)-\varphi(z_0)}, \qquad N(z)=\left\lfloor \log_2 \abs{z-z_0}^{-2} \right\rfloor.
\]
There exist constants $C,c>0$ depending only on $d$ such that
\begin{equation}	\label{eq2219thu}
\abs{H_\varphi^\Omega(z)-\varphi(z_0)}\le\inf_{m \ge 1} \left\{\omega_\varphi(z_0,2^{-m/2})+C\norm{\varphi}_\infty \exp\left(-c(\mathscr W_{N(z)}-\mathscr W_{m-1})_+\right)\right\}
\end{equation}
for every $z\in\Omega$, where the infimum is over integers $m \ge 1$.
In particular,
\begin{equation}	\label{eq1446wed}
\abs{H_\varphi^\Omega(z)-\varphi(z_0)}\le C\left\{\omega_\varphi(z_0,e^{-c\mathscr W_{N(z)}})+\norm{\varphi}_\infty e^{-c\mathscr W_{N(z)}}\right\}.
\end{equation}
\end{theorem}

\begin{remark}
If $z_0$ is regular, \eqref{eq1446wed} gives a modulus of continuity at $z_0$ since $\mathscr W_{N(z)}\to\infty$ as $z\to z_0$.
By subadditivity of thermal capacity and a dyadic decomposition, \eqref{eq2219thu} also yields
\[
\abs{H_\varphi^\Omega(z)-\varphi(z_0)} \le C\left\{\omega_\varphi(z_0,r)+\norm{\varphi}_\infty
\exp\left(-c\int_{\abs{z-z_0}}^r \frac{\capacity(\Omega^\complement\cap\overline{\mathcal E(z_0,s^2)})}{s^{d+1}}\,ds\right)\right\}
\]
for every $z\in\Omega$ with $0<\abs{z-z_0}<r \le 1$.

If $\varphi$ vanishes near $z_0$, the rate of boundary decay is governed by the geometry of $\Omega$.
Indeed, \eqref{eq2219thu} gives $\abs{H_\varphi^\Omega(z)} \lesssim e^{-c\mathscr W_{N(z)}}$.
In particular, $H_\varphi^\Omega$ is Dini continuous at $z_0$ if $\sum_{n=1}^\infty e^{-c\mathscr W_n}<\infty$.
It is H\"older continuous at $z_0$ if
\begin{equation}	\label{eq0902fri}
\capacity(\Omega^\complement\cap\overline{\mathcal E(z_0,s^2)})\ge \sigma s^d
\end{equation}
for some $\sigma>0$ and all sufficiently small $s>0$.
More generally, under \eqref{eq0902fri}, $\alpha$-H\"older continuity ($0<\alpha\le1$) of $\varphi$ at $z_0$ implies $\beta$-H\"older continuity of $H_\varphi^\Omega$ at $z_0$ for some $\beta\in(0,\alpha)$, with respect to the parabolic distance.

Condition \eqref{eq0902fri} holds, for instance, if $\abs{\Omega^\complement\cap\mathcal E(z_0,s^2)}\ge\delta s^{d+2}$ for some $\delta>0$ and all sufficiently small $s>0$, since $\capacity(K)\ge s^{-2}\abs{K}$ for compact $K\subset\overline{\mathcal E(z_0,s^2)}$.
\end{remark}

The proof of Theorem~\ref{thm01} combines an explicit heat-ball barrier, a Carleson-type capacitary packing estimate, and a stopping-time argument to obtain an exponential capacitary estimate in terms of the Wiener series.
The first gives a shorter proof of a key lemma in \cite{EG}, while the latter two provide a quantitative refinement of the argument of Evans and Gariepy.
The main estimates are established in Section~\ref{sec2} and used in Section~\ref{sec3} to prove the theorem.

\section{Capacitary estimates}	\label{sec2}
\subsection{Thermal capacity and a heat-ball estimate}
A function $u\in C(\Omega)$ is caloric in $\Omega$ if $Pu=0$ in the sense of distributions.
A function $v:\Omega\to(-\infty,\infty]$ is supercaloric if it is lower semicontinuous, finite on a dense subset of $\Omega$, and satisfies comparison with caloric functions on cylinders compactly contained in $\Omega$, with boundary inequalities imposed on their parabolic boundaries.

For a compact set $K\subset\bR^{d+1}$, define its reduced function by
\[
\red_K=\inf\{v\ge0:v\text{ is supercaloric in }\bR^{d+1},\; v\ge1\text{ on }K\},
\]
and denote its lower semicontinuous regularization (balayage) by $\bal_K$.
Recall from \cite[Lemmas~2.1 and 2.2]{EG} that thermal capacity is monotone, subadditive, finite on compact sets, invariant under time reversal, homogeneous of degree $d$ under parabolic dilations, and vanishes on singletons.
For every compact set $K$, we have
\[
0\le\bal_K\le\red_K\le1,\qquad \red_K=1\text{ on }K,\qquad \bal_K=1\text{ on }K^\circ.
\]
In $K^\complement$, we have $\bal_K=\red_K$ and $P\bal_K=0$.
The equilibrium measure $\mu_K\ge 0$ supported on $K$ satisfies
\begin{equation}	\label{eq0913sat}
\bal_K=\Gamma*\mu_K,\qquad P\bal_K=\mu_K,\qquad \mu_K(K)=\capacity(K),
\end{equation}
where the second identity holds in the sense of distributions.

If compact neighborhoods $K_i$ decrease to $K$, then $\red_{K_i}\downarrow\red_K$, since for every $v$ admissible for $\red_K$ and $\delta>0$, the function $(1+\delta)v$ is admissible for $\red_{K_i}$ for sufficiently large $i$ by lower semicontinuity.

For $r>0$, write $\mathcal E(r)=\mathcal E(0,r)$ and define the truncated heat ball
\[
\mathcal D(r)=\mathcal E(r)\cap\{t>-3r/4\}=\{(t,x):-3r/4<t<0,\ \abs{x}^2<2d(-t)\log(r/(-t))\}.
\]
We give a shorter proof of Lemma~3.2 in \cite{EG} using a different barrier.

\begin{lemma}	\label{lem01}
Let $r>0$ and let $u\ge0$ be caloric in $\mathcal D(2r)$ and continuous on $\overline{\mathcal D(2r)}\setminus\{0\}$.
There exists $\gamma=\gamma(d)>0$ such that
\[
\inf_{\mathcal E(3r/4)}u\ge\gamma\fint_{\abs{x}^2<3dr/4}u(-3r/2,x)\,dx.
\]
\end{lemma}

\begin{proof}
By scaling, we may assume $r=1$. Set
\[
\alpha=\fint_{\abs{x}^2<3d/4}u(-3/2,x)\,dx.
\]
Write $v(s,x)=u(-s,x)$, so that $P^\ast v=0$, where $P^\ast=-\partial_s-\Delta$, and put
\[
g(s,x)=s^{-d/2}e^{-\abs{x}^2/4s},\qquad h=g^{-1/2}=s^{d/4} e^{\abs{x}^2/8s}.
\]
Direct differentiation gives
\[
P^\ast(1-h)=h\left(\frac{d}{2s}-\frac{\abs{x}^2}{16s^2}\right).
\]
If $P^\ast(1-h)>0$, then $\abs{x}^2/s<8d$ and hence $h\le e^d s^{d/4}$. Thus
\[
P^\ast(1-h)\le\tfrac d2 e^d s^{d/4-1}.
\]
Fix $T\in(0,3/4)$, depending only on $d$, so small that
\[
2e^d T^{d/4}\le(1-(3/4)^{d/4})/2,
\]
and set
\[
w(s,x)=1-h(s,x)+2e^d(s^{d/4}-T^{d/4}).
\]
Then $P^\ast w\le0$, $w\le0$ on $\{0<s\le T,\ g=1\}$, and $w\le1$ on $\{s=T,\ g\ge1\}$.
Moreover,
\begin{equation}	\label{eq1502wed}
w \ge \tfrac{1}{2}(1-(3/4)^{d/4})\quad\text{on }\;\{0<s\le T,\ g\ge(4/3)^{d/2}\}.
\end{equation}

Since $3d/4<3d\log(4/3)$, fix $R>\sqrt{3d/4}$ and $t_1\in(-3/2,-3/4)$, depending only on $d$, such that
\[
(-3/2,t_1]\times\overline B_R(0) \subset\mathcal D(2).
\]
Let $G=G(t,x,y)$ denote the Dirichlet heat kernel in $B_R(0)$.
By comparison and positivity of $G$,
\[
u(t_1,0) \ge \int_{\abs{y}^2<3d/4} G(t_1+3/2,0,y)u(-3/2,y)\,dy \ge c\alpha,
\]
where $c=c(d)>0$ since $R$ and $t_1$ depend only on $d$.
The sets $\{-T\}\times\{x:g(T,x)\ge1\}$ and $\overline{\mathcal E(3/4)}\cap\{t\le-T\}$ are compactly contained in $\mathcal D(2)\cap\{t>t_1\}$.

Since its spatial radius is concave in time, $\mathcal D(2)$ is convex.
A Harnack chain argument starting from $(t_1,0)$ yields
\begin{equation} \label{eq1520wed}
\begin{aligned}
u(-T,x)&\ge c\alpha &&\text{for }\; g(T,x)\ge 1,\\
u&\ge c\alpha &&\text{on }\;\overline{\mathcal E(3/4)}\cap\{t\le-T\},
\end{aligned}
\end{equation}
where $c=c(d)>0$.
Since $v(T,x)=u(-T,x)$, the maximum principle for $P^\ast$ on $\{\epsilon<s<T,\ g>1\}$ yields $v\ge c\alpha w$ in this region.
Letting $\epsilon\downarrow 0$ and combining \eqref{eq1502wed} and \eqref{eq1520wed}, we obtain the assertion.
\end{proof}

\subsection{A capacitary packing estimate}
Write $\mathcal A_j=\mathcal A(0,2^{-j})$ for $j\in\mathbb Z$.
Note that
\[
\mathcal A_j\setminus\{0\}\subset\mathcal E(3\cdot 2^{-j-1}).
\]
For a constant $\lambda>\max\{\sqrt{2},\sqrt{3d/2}\}$ to be fixed below, set
\begin{equation}	\label{eq1047thu}
\mathcal S_j=(-3\cdot 2^{-j},3\cdot 2^{-j}) \times B_{\lambda2^{-j/2}}(0),\qquad j\in\mathbb Z.
\end{equation}

\begin{lemma}	\label{lem02}
There exists a $\lambda=\lambda(d)>\max\{\sqrt{2},\sqrt{3d/2}\}$ such that the following holds.
Let $m\le n$ be integers and let $E_j\subset\mathcal A_j\setminus\{0\}$, $m\le j\le n$, be compact sets.
Let $v\in C^\infty(\bR^{d+1})$ be nonnegative and satisfy $Pv\le0$. Suppose that $v=0$ on $\bigcup_{j=m}^n E_j$.
Set
\[
\alpha_j=\fint_{\abs{x}^2<3d\,2^{-j-1}} v(-3\cdot2^{-j},x)\,dx,\qquad M_j=\sup_{\mathcal S_j}v,\qquad j\in\mathbb Z.
\]
Then
\begin{equation}	\label{eq1521wed}
\sum_{j=m}^n\alpha_j2^{jd/2}\capacity(E_j)\le CM_{m-1},
\end{equation}
where $C=C(d)$.
\end{lemma}

\begin{proof}
Set $f=-Pv\ge0$ and define
\[
w_j=\Gamma*(\mathbbm{1}_{\mathcal E(2^{2-j})}f),\qquad h_j=v+w_j.
\]
Then $h_j\ge v$, $h_j$ is continuous, and $h_j$ is caloric in $\mathcal E(2^{2-j})$.
Applying Lemma~\ref{lem01} to $h_j$ with $r=2^{1-j}$ and using $h_j\ge v$, we obtain $h_j\ge\gamma\alpha_j$ in $\mathcal E(3\cdot2^{-j-1})$.
Since $E_j\subset\mathcal E(3\cdot2^{-j-1})$ and $v=0$ on $E_j$, we have $w_j\ge\gamma\alpha_j$ on $E_j$.

Let $\sigma$ be a nonnegative Radon measure supported on $E_j$ with 
\[
\int_{\bR^{d+1}} \Gamma(z-\zeta)\,d\sigma(z)\le 1,\qquad  \zeta\in\bR^{d+1}.
\]
By Tonelli's theorem,
\[
\gamma\alpha_j\sigma(E_j)\le\int_{E_j}w_j\,d\sigma
=\int_{\mathcal E(2^{2-j})}f(\zeta)\int_{E_j}\Gamma(z-\zeta)\,d\sigma(z)\,d\zeta \le\int_{\mathcal E(2^{2-j})} f.
\]
The supremum of $\sigma(E_j)$ over such measures is the thermal capacity of the time reflection of $E_j$, which equals $\capacity(E_j)$ by the invariance of thermal capacity under time reversal.
Hence
\[
\alpha_j2^{jd/2}\capacity(E_j) \le C2^{jd/2}\int_{\mathcal E(2^{2-j})} f.
\]
Since
$\mathcal E(2^{2-j})=\{\zeta:\Gamma(-\zeta)>(16\pi)^{-d/2}2^{jd/2}\}$,
we have the pointwise estimate
\[
\sum_{j=m}^n 2^{jd/2} \mathbbm{1}_{\mathcal E(2^{2-j})}(\zeta) \le \frac{(16\pi)^{d/2}}{1-2^{-d/2}} \Gamma(-\zeta) \mathbbm{1}_{\mathcal E(2^{2-m})}(\zeta).
\]
Upon integration against $f$, we obtain
\[
\sum_{j=m}^n \alpha_j2^{jd/2}\capacity(E_j) \le Cw_m(0).
\]

To estimate $w_m(0)$, note that
\[
w_m(0)\ge(16\pi)^{-d/2}2^{md/2}\int_{\mathcal E(2^{2-m})} f.
\]
Consider the cylinder
\[
Q=(-3\cdot 2^{1-m},0)\times B_{\lambda2^{-(m-1)/2}}(0)
\]
and denote its parabolic boundary by $\partial_p Q$.
The function $w_m$ vanishes on the bottom of $Q$.
Since the spatial radius of $\mathcal E(2^{2-m})$ is at most $\sqrt{8d/e}\,2^{-m/2}$ and $\lambda>\sqrt{3d/2}$, the bound $\Gamma(t,x)\le C(d)\abs{x}^{-d}$ for $x\ne 0$ gives on the lateral boundary of $Q$
\[
w_m\le C\lambda^{-d}2^{md/2} \int_{\mathcal E(2^{2-m})}f \le C\lambda^{-d}w_m(0),
\]
where $C=C(d)$.
Fix $\lambda=\lambda(d)>\max\{\sqrt{2},\sqrt{3d/2}\}$ sufficiently large that $C\lambda^{-d}\le 1/2$.
Since
\[
P(v+w_m)=-\mathbbm{1}_{\mathcal E(2^{2-m})^\complement}f\le 0
\]
in $Q$ in the sense of distributions and $v\le M_{m-1}$ on $\partial_p Q$ by continuity, the maximum principle gives
\[
v(0)+w_m(0)\le M_{m-1}+\tfrac12 w_m(0).
\]
Thus $w_m(0) \le 2M_{m-1}$, and \eqref{eq1521wed} follows.
\end{proof}

\subsection{A stopping-time argument}
\begin{lemma}	\label{lem03}
Let $\mathcal S_j$ be as in \eqref{eq1047thu}.
There exist $C=C(d)>0$ and $c=c(d)>0$ such that, for every compact $K\subset\bR^{d+1}$ and integers $m\le n$,
\begin{equation}	\label{eq1518wed}
\sup_{\mathcal S_{n+1}}\,(1-\bal_K) \le C\exp\left(-c\sum_{j=m}^n 2^{jd/2}\capacity(K\cap\mathcal A_j)\right).
\end{equation}
\end{lemma}

\begin{proof}
As in \cite[pp.~311--312]{EG}, we first establish a contraction inequality for the supremum over nested cylinders.
Let $E_j$ and $v$ be as in Lemma~\ref{lem02}, with the notation therein.
There exists $\eta=\eta(d)\in(0,1)$ such that
\begin{equation}	\label{eq1516wed}
M_{j+1}\le\eta M_j+(1-\eta)\alpha_j,\qquad j\in\mathbb Z.
\end{equation}
Indeed, $M_j-v$ is nonnegative on $\overline{\mathcal S_j}$ and supercaloric in $\mathcal S_j$.
Let $G=G(t,x,y)$ denote the Dirichlet heat kernel in $B_\lambda(0)$.
By parabolic scaling,
\[
G_j(t,x,y)=2^{jd/2}G(2^jt,2^{j/2}x,2^{j/2}y)
\]
is the Dirichlet heat kernel in $B_{\lambda2^{-j/2}}(0)$.
By comparison, for $(t,x)\in\mathcal S_{j+1}$,
\[
M_j-v(t,x)\ge \int_{\abs{y}^2<3d\,2^{-j-1}} G_j(t+3\cdot2^{-j},x,y)(M_j-v(-3\cdot2^{-j},y))\,dy.
\]
For $(t,x)\in\mathcal S_{j+1}$ and $\abs{y}^2<3d\,2^{-j-1}$, the arguments of $G(2^jt+3,2^{j/2}x,2^{j/2}y)$ lie in a fixed
compact subset of $(0,\infty)\times B_\lambda(0)\times B_\lambda(0)$.
Thus, the positivity of $G$ and the scaling identity give
\[
G_j(t+3\cdot2^{-j},x,y)\ge c2^{jd/2},
\]
where $c=c(d)>0$.
Reducing $c$ if necessary, we obtain
\[
M_j-v\ge(1-\eta)(M_j-\alpha_j)\qquad\text{in }\;\mathcal S_{j+1},
\]
which proves \eqref{eq1516wed}.

Set \[\kappa_j=2^{jd/2}\capacity(E_j)\] and fix $A_\ast=A_\ast(d)\ge 2C$, where $C$ is the constant in \eqref{eq1521wed}.
If $\sum_{j=m}^n \kappa_j\ge A_\ast$, then by Lemma~\ref{lem02} there exists $j_0\in\{m,\ldots,n\}$ such that
\[
\alpha_{j_0}\le CM_{m-1}/A_\ast\le M_{m-1}/2.
\]
Hence, by monotonicity of $M_j$ and \eqref{eq1516wed} with $j=j_0$,
\begin{equation}	\label{eq1503wed}
\sum_{j=m}^n \kappa_j\ge A_\ast \implies M_{n+1}\le\theta M_{m-1},\qquad \theta=(1+\eta)/2<1.
\end{equation}
By scaling and monotonicity of capacity, $\kappa_j\le A_0$ for some $A_0=A_0(d)$.

Set $m_1=m$ and define recursively
\[
n_i=\min\left\{\ell\in\{m_i,\ldots,n\}:\sum_{j=m_i}^\ell \kappa_j\ge A_\ast\right\},\qquad m_{i+1}=n_i+2,
\]
as long as the defining set is nonempty.
By the minimality of $n_i$ and the bound $\kappa_j\le A_0$,
\[
A_\ast\le\sum_{j=m_i}^{n_i}\kappa_j< A_\ast+A_0.
\]
Let $N$ be the number of blocks constructed above (possibly $N=0$).
Since there are at most $N$ skipped indices and the sum of the remaining terms is less than $A_\ast$, we obtain
\begin{equation}	\label{eq1712fri}
\sum_{j=m}^n \kappa_j\le N(A_\ast+2A_0)+A_\ast.
\end{equation}
Since $m_{i+1}-1=n_i+1$ and $\{M_j\}$ is nonincreasing, iteration of \eqref{eq1503wed} yields
\[
M_{n+1}\le\theta^N M_{m-1}.
\]
By \eqref{eq1712fri} and the definition of $\kappa_j$, we obtain
\begin{equation}	\label{eq1459wed}
M_{n+1}\le C M_{m-1}\exp\left(-c\sum_{j=m}^n 2^{jd/2}\capacity(E_j)\right).
\end{equation}

Now choose compact $E_j\subset (K\cap\mathcal A_j)\setminus\{0\}$ and let $\tilde K$ be a compact neighborhood of $K$.
For a standard nonnegative mollifier $\rho_\epsilon$, set 
\[
v_\epsilon=(1-\bal_{\tilde K})*\rho_\epsilon.
\]
Then $0\le v_\epsilon\le1$, $Pv_\epsilon=-\mu_{\tilde K}*\rho_\epsilon\le 0$, and $v_\epsilon=0$ on $K$ for sufficiently small $\epsilon$, since $\bal_{\tilde K}=1$ in $\tilde K^\circ$.
By \eqref{eq1459wed},
\[
v_\epsilon\le C\exp\left(-c\sum_{j=m}^n 2^{jd/2}\capacity(E_j)\right)\quad\text{in }\;\mathcal S_{n+1}.
\]
Since $\bal_{\tilde K}$ is smooth and agrees with $\red_{\tilde K}$ in $\tilde K^\complement$, letting $\epsilon\downarrow 0$ gives the same bound for $1-\red_{\tilde K}$ in $\mathcal S_{n+1}\cap\tilde K^\complement$.
This bound also holds on $\mathcal S_{n+1}\cap\tilde K$, since $\red_{\tilde K}=1$ on $\tilde K$.
Letting $\tilde K\downarrow K$ yields the bound for $1-\red_K$ and hence, by lower semicontinuous regularization on the open set $\mathcal S_{n+1}$, for $1-\bal_K$.
Finally, as in \cite[p.~307]{EG}, we may let $\capacity(E_j)\uparrow\capacity(K\cap\mathcal A_j)$ by subadditivity, continuity of thermal capacity under decreasing compact sets, and the fact that singletons have zero capacity.
This gives \eqref{eq1518wed}.
\end{proof}

\section{Proof of the main theorem}	\label{sec3}
Let $\Omega$ and $z_0=(t_0,x_0)$ be as in Theorem~\ref{thm01}.
For $r>0$ let
\[
\mathcal C_r(z_0)=[t_0-r^2,t_0]\times\overline B_r(x_0),\qquad K_m=\Omega^\complement\cap\mathcal C_{2^{-m/2}}(z_0)\quad(m\ge 1).
\]

\begin{lemma}
There exist positive constants $C=C(d)$ and $c=c(d)$ such that
\begin{equation}	\label{eq1500wed}
1-\bal_{K_m}(z) \le C\exp\left(-c(\mathscr W_{N(z)}-\mathscr W_{m-1})_+\right),\qquad m\ge1,\ z\ne z_0.
\end{equation}
\end{lemma}

\begin{proof}
By translation invariance, we may assume $z_0=0$.
By scaling of $\mathcal E(r)$, there exists an integer $h_0=h_0(d)\ge1$ such that $\Omega^\complement\cap\mathcal A_j\subset K_m$ for $j\ge m+h_0$.
Moreover, by scaling and monotonicity of capacity, $2^{jd/2}\capacity(\mathcal A_j)\le C(d)$.
Note that $z\in\mathcal S_{N(z)+1}$ since $\abs{z}^2\le2^{-N(z)}$ and $\lambda>\sqrt{2}$.
If $N=N(z)\ge m+h_0$, Lemma~\ref{lem03} gives \eqref{eq1500wed}, since
\[
\begin{aligned}
\sum_{j=m+h_0}^N2^{jd/2}\capacity(K_m\cap\mathcal A_j)&=\mathscr W_N-\mathscr W_{m+h_0-1}\\
&\ge\mathscr W_N-\mathscr W_{m-1}-C(d)h_0.
\end{aligned}
\]
The case $N<m+h_0$ follows from $1-\bal_{K_m}\le1$ and $(\mathscr W_N-\mathscr W_{m-1})_+ \le C(d)h_0$.
\end{proof}
\begin{proof}[Proof of Theorem~\ref{thm01}]
By translation invariance, we may assume $z_0=0$.
For necessity, we follow Lanconelli's argument \cite{Lanc73} as presented in \cite[Section~4]{EG}.
Suppose that the series in \eqref{eq1501wed} converges.
Let $F_r=\Omega^\complement\cap\mathcal C_r(0)$.
Splitting the integral representation of $\bal_{F_r}(0)$ in \eqref{eq0913sat} over $\mathcal E(2^{-j_0})$ and its complement, we obtain
\[
\bal_{F_r}(0)\le C2^{j_0d/2}r^d+C\sum_{j=j_0}^\infty 2^{jd/2}\capacity(\Omega^\complement\cap\mathcal A_j).
\]
Here we have used $\mu_{F_r}(F_r)\le\capacity(\mathcal C_r(0))=Cr^d$ by scaling and $\mu_{F_r}(\mathcal A_j)\le
\capacity(\Omega^\complement\cap\mathcal A_j)$ by admissibility of the restricted measure $\mu_{F_r}|_{\mathcal A_j}$.
Choosing $j_0$ sufficiently large and then $r$ sufficiently small gives $\bal_{F_r}(0)<1$.
It follows from \cite[Lemma~1.3]{Lanc73} (see also \cite[Lemma~2.3]{EG}) that $z_0=0$ is irregular.

To prove \eqref{eq2219thu}, let $m\ge 1$ and note that for $h\ge 4$,
\[
\abs{z-\zeta}\ge2^{-m/2}/4 \quad\text{if }\;\abs{z}\ge2^{-m/2}/2 \text{ and }\zeta\in K_{m+h}.
\]
By the equilibrium representation \eqref{eq0913sat}, the bound $\Gamma(z-\zeta)\le C\abs{z-\zeta}^{-d}$ for $z\ne\zeta$, and the estimate $\capacity(K_{m+h})\le C2^{-(m+h)d/2}$, we have
\[
\bal_{K_{m+h}}(z)\le C2^{-hd/2}\quad\text{for }\;\abs{z}\ge2^{-m/2}/2.
\]
Fix an integer $h=h(d) \ge 4$ such that
\[
\bal_{K_{m+h}}\le\tfrac{1}{2} \quad\text{on }\;\{z:\abs{z}\ge2^{-m/2}/2\}.
\]
Set
\[
v_m=\omega_\varphi(0,2^{-m/2})+4\norm{\varphi}_\infty(1-\bal_{K_{m+h}}).
\]
Then $v_m \ge 0$ and $P v_m=0$ in $\Omega$.
By the definition of $\omega_\varphi$ and the choice of $h$,
\[
\liminf_{\Omega\ni z\to\zeta}v_m(z)\ge\abs{\varphi(\zeta)-\varphi(0)},\qquad \zeta\in\partial\Omega.
\]
By the Perron method, $\abs{H_\varphi^\Omega-\varphi(0)}\le v_m$ in $\Omega$.
Applying \eqref{eq1500wed} and absorbing the contribution of the fixed shift $h$ into the constant $C$, we obtain \eqref{eq2219thu}.
Since $\mathscr W_{m-1}\le C(d)(m-1)$,
\[
\exp\left(-c(\mathscr W_{N(z)}-\mathscr W_{m-1})_+\right) \le \exp(C_1(m-1)-c\mathscr W_{N(z)}),\qquad C_1=C_1(d)>0.
\]
Choosing $m=\max\{1,\lceil c\mathscr W_{N(z)}/(2C_1)\rceil\}$, we obtain
\[
2^{-m/2}\le e^{-\kappa\mathscr W_{N(z)}},\qquad C_1(m-1)-c\mathscr W_{N(z)}\le-\tfrac{c}{2}\mathscr W_{N(z)},
\]
where $\kappa=c\log 2/(4C_1)$.
This proves \eqref{eq1446wed} with $c$ replaced by $\min\{\kappa,c/2\}$.

If \eqref{eq1501wed} holds, then $\mathscr W_{N(z)}\to\infty$ as $z\to0$ since $N(z)\to\infty$.
Thus the right-hand side of \eqref{eq1446wed} tends to zero, and hence $z_0=0$ is regular.
\end{proof}
\enlargethispage{\baselineskip}

\end{document}